\documentclass[12pt]{amsart}
\usepackage[english,german]{babel}
\usepackage{amsmath}
\usepackage{amssymb}
\usepackage{amsfonts}
\usepackage{cmap}

\begin{document}

\thispagestyle{empty}

\title[On applications of the model spaces\dots]{On applications of the model spaces \\ to the construction of cocyclic perturbations \\ of the semigroup of shifts on the semiaxis}

\author{G.G. Amosov, A.D. Baranov, V.V. Kapustin}

\address{Grigory Amosov,
\newline\hphantom{iii} Steklov Mathematical Institute of RAS,% ћесто работы
\newline\hphantom{iii} Gubkin Str., 8, 
\newline\hphantom{iii} 119991, Moscow, Russia}
\email{gramos@mi.ras.ru}

\address{Anton Baranov,
\newline\hphantom{iii} St. Petersburg State University,% ћесто работы
\newline\hphantom{iii} Staryi Petergof, Universitetskii pr., 28, % јдрес (улица, дом, строение и т.п.)
\newline\hphantom{iii} 198504, St. Petersburg, Russia}%  јдрес (почтовый индекс, город, страна)
\email{anton.d.baranov@gmail.com}% ¬аш электронный адрес дл€ переписки

\address{Vladimir Kapustin,
\newline\hphantom{iii} St. Petersburg Department\\ Steklov Mathematical Institute of RAS,% ћесто работы
\newline\hphantom{iii} Fontanka levee, 27, % јдрес (улица, дом, строение и т.п.)
\newline\hphantom{iii} 191023, St. Petersburg, Russia}%  јдрес (почтовый индекс, город, страна)
\email{kapustin@pdmi.ras.ru}% ¬аш электронный адрес дл€ переписки

%\thanks{\sc G.G. Amosov, A.D. Baranov, V.V. Kapustin, %  ‘.».ќ. авторов на английском €зыке
%On applications of the model spaces to the construction of cocyclic perturbations of the semigroup
%of shifts on the semiaxis.}% название статьи на английском €зыке
\thanks{\copyright \ G.G. Amosov, A.D. Baranov, V.V. Kapustin 2012}
\thanks{\rm The work is supported by the program of RAS "Mathematical
basics of management"\ and the grant of Russian Foundation for Basic Research 11-01-00584-a.}
\thanks{\it Submitted on 20 December 2011}

\sloppy

\maketitle {\small
\begin{quote}
\noindent{\bf Abstract. } We describe a construction of cocyclic perturbations of the
semigroup of shifts on the semiaxis by means of the theory of model spaces. It is shown
that one can choose an inner function that determines the model space so that the
elements of the perturbed semigroup have a prescribed spectral type and differ from the
elements of the initial semigroup by operators from the Schatten--von Neumann class
$\mathfrak{S}_p$, $p>1$. The case of the trace class $\mathfrak{S}_1$ perturbations is
considered separately.

\medskip

\end{quote}
 }

\section{Introduction}%ќсновной текст статьи

Let us assume that $(S_t,\ t\ge 0),$ and $(\tilde S_t,\ t\in {\mathbb R})$ are the
semigroup of shifts in the space  $H=L^2({\mathbb R}_+),$ and the group of shifts (its
unitary dilation) in the space $\tilde H=L^2({\mathbb R})$, defined by formulae
$$
(S_tf)(x)=
\begin{cases} f(x-t), & x>t, \\
              0,  &     0\le x\le t,
\end{cases} \qquad  f\in H,
$$
and
$$
(\tilde S_tg)(x)=g(x-t),  \qquad g\in \tilde H,
$$
respectively. Sometimes it is convenient to assume that the multiplicative group of
the algebra $B(H)$ of bounded operators in the space $H$  is embedded in the
multiplicative group of the algebra  $B(\tilde H)$ in such a way that elements $B(H)$
act on functions $f\in \tilde H$ with a support on the negative semiaxis as the identical
mapping. In this case operators, acting in the space $H$, will be considered   as
operators in $\tilde H$ as well. A strongly continuous family of unitary operators
$(W_t,\ t\ge 0)$ in the space $H$ is called a {\it cocycle} of the semigroup of shifts
$(S_t,\ t\ge 0)$  if the condition
\begin {equation}\label {coc}
W_{t+s}=W_t\tilde S_tW_s\tilde S_{-t},\quad t,s\ge 0,\qquad W_0=I
\end{equation}
 holds true (see \cite {Arv}).
It follows from the condition (\ref {coc}) that the family of isometric operators
$(V_t=W_tS_t,\ t\ge 0)$ in the space $H$ forms a semigroup (i.e. $V_{t+s}=V_tV_s$,
$t,s\ge 0$), which will be called a {\it cocyclic perturbation} of the semigroup of shifts
$(S_t,\ t\ge 0)$.
\smallskip

In the present paper we will show that any cocyclic
perturbation of the semigroup $(S_t)$ is unitarily equivalent to the orthogonal sum
\begin {equation}\label {res}
(V_t)\cong (U_t\oplus S_t),%\quad t\ge 0,
\end{equation}
where $(U_t,\ t\ge 0)$ is a semigroup of unitary operators, and the following two
theorems hold true. Here and in what follows all  semigroups under consideration are
supposed to be strongly continuous; the symbol $\mathfrak {S}_p$ denotes the
Schatten-von Neumann operator ideals.
\medskip

{\bf Theorem 1.} {\it For any semigroup of unitary
operators $(U_t,\ t\ge 0)$ possessing a spectral measure 
which is singular with respect to the Lebesgue measure,
there is a cocycle $(W_t,\ t\ge 0)$, satisfying the condition
$$
W_t-I\in \mathfrak {S}_p
$$
for all $p>1$, where~$(\ref {res})$ holds true for
the cocyclic perturbations $(V_t=W_tS_t,\ t\ge 0)$, and}
\begin{equation}\label{utv}
V_t-S_t\in \mathfrak {S}_1,\qquad t\ge 0.
\end{equation}

As a corollary of Theorem 1 we obtain analogous results for an arbitrary
(not necessarily singular) spectral measure.
\medskip

{\bf Theorem 2.} {\it For any semigroup of unitary
operators $(U_t,\ t\ge 0)$ and for any $p>1$ there is a
cocycle $(W_t,\ t\ge 0)$, satisfying the condition
$$
W_t-I\in \mathfrak {S}_p
$$
for all $p>1$, when the relation~$(\ref {res})$ holds true for the cocyclic perturbation}
$(V_t=W_tS_t,\ t\ge 0)$.

\medskip

In what follows it will be shown (Proposition~10)
that  the condition
$W_t-I\in \mathfrak {S}_1$
never holds true in the considered model of cocyclic perturbations. 
Thus, the results of the work are in a sense optimal.
It is natural to suppose, that this fact holds true in the general case.
\medskip

{\bf Conjecture.} For any cocycle $(W_t,\ t\ge 0),$ 
such that $W_t-I\in\mathfrak S_1$
for all $t\ge 0$, the perturbed semigroup $(V_t=W_tS_t,\ t\ge 0)$ 
is unitarily equivalent to the initial one: $(V_t)\cong (S_t)$.

\medskip

We should note that the problem of the Markov cocyclic perturbations of the 
group of unitary operators connected with the questions under 
consideration is posed in \cite {Amo}, 
and the Markov cocycles possessing the property
$W_t-I\in \mathfrak {S}_2,\ t\ge 0$ are considered in
\cite {Bar1, Bar2}. The property~(\ref {utv}) was considered 
in the article~\cite{ABK}, where perturbations
$(V_t,\ t\ge 0)$ of semigroup shifts $(S_t,\ t\ge 0)$ such  that
$V_t-S_t\in \mathfrak {S}_p,\ p\ge 1$, were investigated. 
A distinguishing feature of the present paper is that the considered 
perturbations  possess additional cocyclic properties demanding consideration 
of unitary dilations of semigroups. The technique applied here 
is analogous to the one in paper \cite {ABK}.

\section {Cocyclic perturbations of the general form}

For any strongly continuous semigroup of isometric operators $(V_t,\ t\ge 0)$
in the Hilbert space $H,$ the  
Wold--Kolmogorov decomposition is defined as follows:
$$
H=H_0\oplus H_1,
$$
\begin{equation}
   \label{Wold}
V_t=U_t\oplus R_t,\qquad t\ge 0,
\end{equation}
where $(U_t,\ t\ge 0)$ is a semigroup of unitary operators in $H_0$, and $(R_t,\ t\ge 0)$ is a semigroup
of completely nonunitary isometric operators in $H_1$, i.e. lacking
nontrivial  invariant
subspaces, where they act  as unitary operators.
\medskip

{\bf Proposition 3.} {\it Let the semigroup of
isometric operators $(V_t,\ t\ge 0)$ be a cocyclic perturbation of
the semigroup of shifts $(S_t,\ t\ge 0)$.
Then  the completely nonunitary part $(R_t,\ t\ge 0)$ in the Wold--Kolmogorov decomposition
$(\ref {Wold})$ is unitarily equivalent to the semigroup of shifts $(S_t,\ t\ge 0)$.}

\medskip

{\bf Remark.} This statement holds true for an arbitrary  semigroup (not necessarily
being a cocyclic perturbation) of isometric operators $(V_t,\ t\ge 0)$,
if we require that
$V_t-S_t\in \mathfrak {S}_p$, $p\ge 1$ (see~\cite {ABK}).

\medskip

{\bf Proof.} Let us define elements $\xi _t\in H,\ t\ge 0$, by the formula
$$
\xi _t(x)= \begin{cases} 1, &  0\le x\le t,  \\
                         0, &   x>t.
\end{cases}
$$
Note, that the family $(\xi _t,\ t\ge 0)$ satisfies
the so-called condition of an {\it additive cocycle} of the semigroup $(S_t,\ t\ge 0)$,
i.e.
$$
\xi _{t+s}=\xi _t+S_t\xi _s,\qquad s,t\ge 0,
$$
and the functions $\xi _{t_1}-\xi _{s_1}$ and $\xi _{t_2}-\xi _{s_2}$ are orthogonal, if
$(s_1,t_1)\cap (s_2,t_2)=\emptyset$. Moreover, linear combinations of elements
$(\xi _s,\ 0\le s\le t)$ generate ${\rm Ker}\, S_t^*$.
Assume that $\tilde \xi _t=W_t\xi _t,\ t\ge 0$.
To prove Proposition 3 it is sufficient to verify that
for the cocyclic perturbation $(V_t=W_tS_t,\ t\ge 0)$
the family of elements $\tilde \xi _t$ has the following properties:

(i) $\tilde \xi _{t+s}=\tilde \xi _t+V_t\tilde \xi _s,\ s,t\ge 0,$

(ii) $\tilde \xi _{t_1}-\tilde\xi _{s_1}$ and $\tilde \xi _{t_2}-\tilde \xi _{s_2}$ are orthogonal,
if $(s_1,t_1)\cap (s_2,t_2)=\emptyset$,

(iii) linear combinations $(\tilde\xi _s,\ 0\le s\le t)$ generate
${\rm Ker}\, V_t^*$.

\noindent Indeed, in this case the restriction  of the semigroup $(V_t,\ t\ge 0)$ to the subspace $H_0$,
generated by ${\rm Ker}V_t^*,\ t\ge 0,$ is unitarily equivalent $(S_t,\ t\ge 0)$,
but the restriction $V_t|_{H_0^{\perp}}$ will be a unitary operator,
because ${\rm Ker}\,V_t|_{H_0^{\perp}}=\{0\},\ t\ge 0$.

We have
\begin{equation}\label{0}
\tilde \xi _{t+s}=W_{t+s}\xi _{t+s}=
W_t\tilde S_tW_s\tilde S_{-t}\xi _t+W_t\tilde S_t(W_s)\tilde S_{-t}S_t\xi _s.
\end{equation}
Note, that
\begin{equation}\label{1}
\tilde S_tW_s\tilde S_{-t}\xi _s=\xi _s,
\end{equation}
whereas $W_sf=f$ for the function with the support 
${\rm supp}\,f\subset {\mathbb R}_-$. On the other hand,
\begin{equation}\label{2}
\tilde S_{t}W_s\tilde S_{-t}S_t\xi _s=\tilde S_tW_s\xi _s=S_t\tilde \xi _s.
\end{equation}
Substituting the relations (\ref {1}) and (\ref {2}) into the equality
(\ref {0}), we obtain the property (i).

Next,
\begin{equation}\label{3}
W_{t+s}\xi _t=W_t\tilde S_t(W_s)\tilde S_{-t}\xi _t=W_t\xi _t,\qquad s,t\ge 0,
\end{equation}
according to (\ref {1}).
Let $\tilde t= \max (t_1,t_2)$; then, taking into account (\ref {3}),
we obtain
$$
\begin{aligned}
 (\tilde \xi _{t_1}-\tilde\xi _{s_1}, \tilde \xi _{t_2}-\tilde \xi _{s_2})
& = (W_{t_1}\xi _{t_1}-W_{s_1}\xi _{s_1},W_{t_2}\xi _{t_2}-W_{s_2}\xi _{s_2}) \\
& = (W_{\tilde t}\xi _{t_1}-W_{\tilde t}\xi _{s_1},W_{\tilde t}\xi _{t_2}-W_{\tilde t}\xi _{s_2})
=(\xi _{t_1}-\xi _{s_1},\xi _{t_2}-\xi _{s_2})=0,
\end{aligned}
$$
if $(s_1,t_1)\cap (s_2,t_2)=\emptyset$. Thus, the property
(ii) is established as well.

Finally, let us consider the equation
\begin{equation}\label{4}
V_t^*f=S_t^*W_t^*f=0.
\end{equation}
It follows from (\ref {4}) that ${\rm supp}\,  W_t^*f\subset [0,t]$.
Hence, $f$ belongs to the closure of the linear
envelope  of elements $(W_t\xi _s,\ 0\le s\le t)$. Since $s\le t$, we have
$W_t\xi _s=W_s\xi _s=\tilde \xi _s$ for such elements
by virtue of the relations (\ref {3}). This completes the proof 
of the property (iii) and of the proposition.
\qed
\medskip

The next property is necessary for the construction of a model of cocycles.
\medskip

{\bf Proposition 4.}
{\it Let $(V_t,\ t\ge 0)$ be the cocyclic perturbation of the  semigroup of shifts $(S_t,\ t\ge 0)$
by the cocycle $(W_t,\ t\ge 0)$. Then, defining the family of the unitary operators
$(W_{-t},\ t\ge 0)$ in the space $\tilde H$ by the formula
\begin{equation}\label {cont}
W_{-t}=\tilde S_{-t}W_t^*S_t,\qquad  t \ge 0,
\end{equation}
we obtain that the family of operators $(\tilde V_t,\ t\in {\mathbb R})$, where
$$
\tilde V_t=W_t\tilde S_t,
$$
generates a group of unitary operators in the space $\tilde H$, and}
$$
\tilde V_tf=
\begin{cases} V_tf, & \quad {\rm supp}\,  f \subset {\mathbb R}_+,\quad t\ge 0, \\
        \tilde S_tf, & \quad {\rm supp}\, f\subset {\mathbb R}_-,\quad t\le 0.
\end{cases}
$$

\smallskip

{\bf Proof.} As usual, let us assume
that actions of the unitary operators $W_t,\ t\ge 0$, fixed
initially in the space
$H$, are extended by the identical action on $f$ with 
${\rm supp}\,f\subset {\mathbb R}_-$. Then the formula (\ref {cont})
provides the prolongation  of the family $(W_t,\ t\ge 0)$ of unitary operators in $\tilde H$ for negative values of the parameter $t$.
Moreover, the property of the cocycle
$$
W_{t+s}=W_t\tilde S_tW_s\tilde S_{-t},\qquad  s,t\in {\mathbb R}
$$
holds, which follows from the formula
$$
I=W_{-t+t}=W_{-t}\tilde S_{-t}W_t\tilde S_t,\qquad  t \ge 0,
$$
resulting from the definition (\ref {cont}). To complete the proof it should be noted, that
for ${\rm supp}\,f\subset {\mathbb R}_-$,
$$
\tilde V_{-t}f=W_{-t}\tilde S_{-t}f=\tilde S_{-t}W_{t}^*f=\tilde S_{-t}f,
\qquad t\ge 0.
$$
\qed

\section {Model of the cocyclic perturbation based on the cogenerator of the semigroup}

We will need some well-known background from the theory 
of one-parameter semigroups (see \cite {Nagy}). A symmetric
(probably unbounded) operator $A=s-\lim_{t\to 0+} \frac {V_t-I}{it}$ 
is called the generator of a strongly continuous group of the isometric 
operators $(V_t,\ t\ge 0).$
An isometric operator $V=(A-iI)(A+iI)^{-1}$ is called a cogenerator of a semigroup.
For an isometric operator to be a cogenerator of some isometric semigroup
it is necessary and sufficient that the number $1$ should not belong to its point spectrum.
The initial semigroup will consist of unitary operators only if
$A$ is a self-adjoint operator, or, equivalently, when $V$ is a
unitary operator such that the point $1$ does not belong to its point spectrum.
If we introduce the functions
\begin{equation}\label{funct}
\varphi _t(z)=\exp \bigg(t\frac {z+1}{z-1}\bigg),\qquad t\ge 0,
\end{equation}
then the semigroup is recovered by means of the cogenerator $V$  as follows:
$V_t=\varphi _t(V),\ t\ge 0$.
Let us note, that the functions $\varphi_t$ are bounded and analytic 
in the unit disc $\mathbb{D}$.

One can readily show that  the cogenerator of the  semigroup of shift operators
$(S_t,\ t \ge 0)$ in the space $H$ is unitarily equivalent to the operator of the
(unilateral) shift $S$ in the Hardy space $K=H^2({\mathbb D})$,
consisting of analytical
in the circle $\mathbb{D}$ functions $f(z)=\sum_{n=0}^{+\infty}c_nz^n$,
for which $\sum_{n=0}^{+\infty} |c_n|^2=\|f\|^2_{L^2(\mathbb {T})}<+\infty$.
Therefore, the Hardy space in the disc is naturally
embedded in the space $\tilde K=L^2({\mathbb T})$
on the circle $\mathbb{T}$. The operator of the shift $S$ in the Hardy space
is given by the formula
\begin{equation}\label{sh}
  (Sf)(z)=zf(z), \qquad f\in K.
\end{equation}

Analogously, the cogenerator of the  group of shifts
in the space $\tilde H$  is unitarily equivalent to the operator 
of the (bilateral) shift $(\tilde S f)(z) = z f(z)$ in the space $\tilde K$, with
the operator $\tilde S$, apparently, being a unitary dilation of the operator $S$.

Assume, that $E$ is a nontrivial invariant subspace of the shift operator
$S$, that is, $SE\subset E$. Then, according to the Beurling theorem (see \cite {Nik}),
$E=\theta H^2({\mathbb D})$ for some inner function $\theta\in H^{\infty}({\mathbb D})$
(i.e., the function which is analytic and bounded in the unit disc
$\mathbb D$ with nontangential boundary values such that
$|\theta (z)|=1$ almost everywhere on $\mathbb T$). The orthogonal complement
$K_{\theta }=H^{2}({\mathbb D})\ominus \theta H^{2}({\mathbb D})=E^{\perp}$ is usually called a
{\it model} space. The next proposition describes the model of cocyclic perturbation,
applied in the present paper.

\medskip

{\bf Proposition 5.} {\it A cogenerator of any
cocyclic perturbation of the  semigroup of shifts on the semiaxis 
is unitarily equivalent to the isometric operator $V$ in the space
$K=H^2({\mathbb D})$, for which there is an inner function $\theta $,
such that
\begin{equation}
  \label{Wold1}
  V=U\oplus S|_{E},
\end{equation}
where $S|_{E}$ is the restriction of the operator of the shift $S$ to
the invariant space defined by the function $\theta $, and $U$ is
a unitary operator in the model space $K_{\theta}$, which is 
the cogenerator of the unitary part of the
Wold--Kolmogorov decomposition of the cocyclic perturbation.}

\medskip

{\bf Proof.} For the cogenerator of the cocyclic perturbation
$V$ in the space $K,$ one has a determined  Wold--Kolmogorov decomposition
$K=K_0\oplus K_1$ such that
$V|_{K_0}$ is a unitary operator and the restriction $V|_{K_1}$
is a completely  nonunitary isometric operator.
It follows from Proposition 3 that the restriction $V|_{K_1}$ is unitarily equivalent to the shift operator $S$.
Therefore, in our model situation, one can use 
as $V|_{K_1}$ the restriction 
$S|_{E}$ to any invariant subspace $E$, selected so that
the equality ${\rm dim}\,K_{\theta}={\rm dim}\, K_0$ holds 
true for the corresponding model space $K_{\theta} = E^\perp,$ 
which completes the proof. \qed
\medskip

The following statement  results directly  from Proposition 5.
\medskip

{\bf Corollary 6.} {\it The cogenerator of the semigroup of unitary operators $(\tilde V_t,\ t\ge 0)$,
defining the cocycle according to  Proposition~4, is unitary equivalent to the operator $\tilde V$
in the space $\tilde K=L^2({\mathbb T})$, possessing the properties}
$$
\begin{aligned}
\tilde Vf=Vf,\qquad & f\in K=H^2({\mathbb D}), \\
(\tilde V^*f)(z)=\overline zf(z),\qquad
& f\in \tilde K\ominus K=L^2({\mathbb T})\ominus H^2({\mathbb D}).
\end{aligned}
$$

\section {Perturbation model based on the Clark measures}

Let $U$ be the unitary part in the Wold--Kolmogorov decomposition
(\ref {Wold1}) of the cogenerator of the cocyclic perturbation.
In this section we will be interested in the case when $U$ is 
unitarily equivalent to
the operator of multiplication by $z$ in the space $L^2(\mu)$,
%$$
%(Uf)(z)=zf(z),\qquad  f\in L^2(\mu),
%$$
with the measure $\mu$ being singular with  respect to the Lebesgue measure. Note, that
$U$ a cogenerator of the semigroup according to the condition  and  therefore  the number $1$ does not belong to its
point spectrum. Operators of multiplication by $z$ in the spaces
$L^2(\mu)$ and $L^2(\tilde \mu)$ are unitarily equivalent,
if the measures $\tilde \mu$ and $\mu$ are mutually absolutely continuous.
Multiplying the measure $\mu$ by a positive weight, one can make it
satisfy the following auxiliary  condition,
taking  an important part in what follows:
\begin{equation}
 \label{usl}
 \int \limits _{\mathbb T}\frac {d\mu(\xi)}{|1-\xi|^q}<+\infty
\end{equation}
for some $q>3$.

Let $\mu$ be a finite singular Borel measure on the unit circle.
Define the inner function $\theta$ by the formula
\begin{equation}
\label{Clark}
\frac {1+\theta (z)}{1-\theta (z)}=
\int \limits _{\mathbb T}\frac {\xi +z}{\xi -z}d\mu (\xi).
\end{equation}
Then the operator $\Omega$, given on $L^2(\mu)$ by the formula
\begin{equation}
  \label{omega}
  (\Omega f)(z)=(1-\theta (z))\int \limits _{\mathbb T}
  \frac {f(\xi)d\mu (\xi)}{1-\overline {\xi}z},
\end{equation}
is a unitary operator from $L^2(\mu)$ onto $K_{\theta}$.
Moreover, the unitary operator $U$ in $L^2(\mu)$
transforms into the unitary operator $\tilde U$ in the model space $K_{\theta}$
such that
\begin{equation}\label{unit}
\tilde Uf = \Omega U\Omega^* f = zf +(f,g) (1-\theta ),\ f\in K_{\theta},
\end{equation}
where
$$
g(z)=\frac {\theta (z)-\theta (0)}{z(1-\theta (0))} \in K_{\theta},
$$
and therefore the operators $U$ and $\tilde U$ are unitarily equivalent, see \cite {Clark}.

The operator (\ref {unit})
is the restriction to the model space $K_{\theta}$
of the isometric operator $V$, acting in the space $K$ by formula
\begin{equation}\label{mod1}
(Vf)(z)=zf(z)+(f,g)(1-\theta (z)),\qquad  f\in K.
\end{equation}
The unitary dilation of the operator (\ref {mod1}) will be the operator
\begin{equation}\label{mod2}
(\tilde Vf)(z)=zf(z)+(f,g)(1-\theta (z))
-(f,\overline z)(1-\overline{\theta (1)} \theta (z)),\quad f\in \tilde K.
\end{equation}
Note that
$$
(\tilde V^*f)(z)=\overline zf(z),\qquad  f\in \tilde K\ominus K.
$$
Therefore, according to Proposition 5 and Corollary 6,
the following statement is proved.

\medskip

{\bf Proposition 7.} {\it The formulae $(\ref {mod1})$, $(\ref {mod2})$ 
define the model of the cogenerator of a cocyclic perturbation 
in the case, when the unitary part of the cogenerator 
in the Wold--Kolmogorov decomposition is
unitarily equivalent to the operator of multiplication by $z$
in the space $L^2(\mu)$ with the measure $\mu$, which is 
singular with respect to the Lebesgue measure.}

\section {Closeness of cocyclic perturbation}

Let us apply the function (\ref {funct}) to the model
cogenerator $V$ of the semigroup of isometric operators $(V_t,\ t\ge 0)$.
The isometric operator $V$ is the restriction of the unitary operator $\tilde V$
defined by the formula (\ref {mod2})
to the space $K = H^2$. Recall that the symbols $S$ and $\tilde S$
define the shift operators  on $K$ and $\tilde K,$ respectively.
Then the cocycle $(W_t,\ t\ge 0)$ satisfies the equality
$$
\varphi _t(\tilde V)-\varphi _t(\tilde S) = (W_t-I)\tilde S_t, \quad t\ge 0.
$$
Therefore, inclusion of the difference $W_t-I$ into the ideals
$\mathfrak {S}_p$ proves to be equivalent to the corresponding inclusion
for the differences $\varphi _t(\tilde V)-\varphi _t(\tilde S)$.
The properties of the operators
$\varphi _t(\tilde V)-\varphi _t(\tilde S)$ are determined in their turn
by properties of the spectral measure $\mu$ of the unitary operator
$(\ref {unit})$, i.e., by its smallness (smoothness)
at the point 1.

We will need the following statement,
proved in \cite{ABK} (Proposition 7.2).

\medskip

{\bf Proposition 8.} {\it Let the spectral measure of the unitary
operator $(\ref {unit})$ satisfy the condition
\begin{equation}
 \label{uslovie}
 \mathfrak {M}_q(\mu) = \int \limits _{\mathbb T}
 \frac {d\mu(\xi)}{|1-\xi |^q}<+\infty
\end{equation}
for some $q>3$. Then
$$
\varphi _t(V)-\varphi _t(S)\in \mathfrak {S}_1,\quad t\ge 0,
$$
with
$$
\|\varphi _t(V) - \varphi _t(S)\|_{\mathfrak {S}_1} \le C_q t^{1/2}
(\mathfrak {M}_q(\mu))^{1/2},
$$
where the constant $C_q$ depends only on $q$.}

\medskip

The key role in the proof of the Theorems 1 and 2 is played by the following proposition,
allowing one to estimate components of the unitary  dilation.
In this case we are not able to obtain the inclusion of
$\varphi _t(\tilde V)-\varphi _t(\tilde S)\in \mathfrak {S}_1$,
but the difference may belong to the ideals
$\mathfrak {S}_p$ for all $p>1$.

\medskip

{\bf Proposition 9.} {\it Let the spectral measure of the unitary
operator $(\ref {unit})$ satisfy the condition
$(\ref{uslovie})$ for some $q>3$. Then
$$
\varphi _t(\tilde V)-\varphi _t(\tilde S)\in \mathfrak {S}_p,\qquad
p>q' = \frac {q}{q-1}, \quad t\ge 0,
$$
with
$$
\|\varphi _t(\tilde V)-\varphi _t(\tilde S)\|_{\mathfrak {S_p}}\le
\omega (\mathfrak {M}_q(\mu)),
$$
where $\omega$ is a positive function such that
$\omega (r) \to 0$ when $r \searrow 0$.}

\medskip

{\bf Proof.}
The proof of Proposition 9 consists of several steps.
At the first step we will consider components of the operator
$\varphi _t(\tilde V)-\varphi _t(\tilde S)$ with respect  to
some canonical representation of the space $\tilde K$
and will see that all the components, except one,
belong to the ideal $\mathfrak {S}_1$
due  to Proposition 8. Then, we will show that
the remaining  component is unitarily equivalent (after conformal transformation
to the upper half-plane) to the operator of multiplication by a certain
function in the PaleyЦ-Wiener space.
This will allow us  to reduce the problem to   the question
of  describing measures (weights), such that  the embedding operator of  
the PaleyЦ-Wiener space belongs to the ideal $\mathfrak {S}_p$.
To complete the proof we apply a theorem due to O.G. Parfenov
\cite{Par}.
\smallskip
\\
{\it \bf Step 1. Analysis of  components of the unitary dilation.}
Let us consider the matrix of the operator $\varphi_t(\tilde V)-
\varphi_t(\tilde S)$ with respect to  the expansion
$\tilde K= H^2_- \oplus K_\theta \oplus \theta H^2$,
where $H^2_- = L^2(\mathbb{T}) \ominus H^2$.
One can readily see that  
all the components, except one, belong to the class
$\mathfrak{S}_1$. Indeed, 
the statement follows from Proposition 8 for the block
$ K_\theta \oplus \theta H^2 \to  K_\theta\oplus\theta H^2$.
Proceeding to the conjugate operator, we come to the conclusion that
the block $H^2_-\oplus K_\theta\to H^2_-\oplus K_\theta$
is also included into $\mathfrak{S}_1$.
By its construction the component $H^2\to H^2_-$ is equal to zero.
Therefore, we only need to consider the component, corresponding to the operator
$H^2_-\to\theta H^2$. Moreover,  note that both operators
$\varphi_t(\tilde V)$ and $\varphi_t(\tilde S)$
on the space $\bar\varphi_t H^2_-$ act as operators of multiplication
by $\varphi_t$, and, consequently, $\varphi_t(\tilde V)-\varphi_t(\tilde S)=0$
on $\bar\varphi_t H^2_-$.
It remains only to study the action  of the operator $\varphi_t(\tilde V)-\varphi_t(\tilde S)$
on the subspace $\bar\varphi_t H^2\ominus H^2=\bar\varphi_t K_{\varphi_t}$.
Let us denote the restriction of the operator $\varphi_t(\tilde V)-\varphi_t(\tilde S)$
to the subspace $\bar\varphi_t K_{\varphi_t}$ by $Q: \bar\varphi_t K_{\varphi_t} \to H^2$.
\smallskip
\\
{\it \bf Step 2. Inclusion of the component $Q$ into the ideals $\mathfrak{S}_p$.}
Let us show that for $v\in K_{\varphi_t}$ the following equality holds true
\begin{equation}
\label{dil0}
Q(\bar\varphi_t v)=-(1- \overline{\theta (1)} \theta)v.
\end{equation}
If $u\in H^2_-$, then
for the arbitrary function $\varphi\in H^\infty$ there is the equality
\begin{equation}
\label{dil}
P_+\varphi(\tilde V) u = \overline{\theta (1)} \theta \cdot P_+ (\varphi u),
\quad P_-\varphi(\tilde V) u = P_-(\varphi u),
\end{equation}
where the symbols $P_+$ and $P_-$ denote projectors in the space $L^2(\mathbb{T})$
on the subspace $H^2$ and $H^2_-,$ respectively.
Indeed, this equality is easily verified for the case when
$\varphi(z)=z^n$, $n>0$, and $u(z) = z^m$, $m<0$.
Due to its linearity and continuity the equality (\ref{dil}) holds true for
all $u\in H^2_-$ and $\varphi(z)=z^n$, $n>0$.
Finally, due to its linearity and $\ast$-weak continuity,
the equality (\ref{dil}) holds   for the arbitrary function
$\varphi\in H^\infty$ as well.

Since $\varphi_t(\tilde S)u=\varphi_tu$, the equality
(\ref{dil}) entails 
$$
\big(\varphi_t(\tilde V)-\varphi_t(\tilde S)\big)u =
(\overline{\theta (1)} \theta-1) \cdot P_+ (\varphi_t u), \qquad u\in H^2_-.
$$
Substituting $u=\bar\varphi_t v$, we  obtain the equality (\ref{dil0}).

Therefore,
the inclusion of  $\varphi_t(\tilde V)-\varphi_t(\tilde S) \in
\mathfrak {S}_p$ is equivalent to the inclusion
\begin{equation}
  \label{inc}
   M_{1-\overline \theta (1)\theta}
|_{H^2\ominus \varphi _tH^2} \in \mathfrak {S}_p,
\end{equation}
where the symbol $M_g$ denotes the operator of multiplication by the function
$g\in L^\infty(\mathbb{T})$.
\smallskip
\\
{\it \bf Step 3. Transformation  into the half-plane.}
It will be convenient to prove the inclusion (\ref {inc}),
making a ``unitary transformation'' from the unit disc into the upper half-plane
$\mathbb C_+=\{z: {\rm Im}\, z>0\}$.
Let us assume that
$$
\Theta (z)=\theta\bigg(\frac {z-i}{z+i}\bigg).
$$ 
Then $\Theta (z)$ becomes an inner function in the upper half-plane:
$\Theta \in H^{\infty}({\mathbb C}_+)$, and
${|\Theta (x)|=1}$ for almost every $x\in {\mathbb R}$, where 
the values of the function $\Theta$ on the real line are considered 
as nontangential boundary values.
Defining  the measure $\nu $ on the real line by the condition
$$
d\mu (\xi)=\frac {d\nu (x)}{\pi (1+x^2)},\qquad \xi =\frac {x-i}{x+i},
$$
we  obtain
$$
\frac {1-\Theta (z)}{1+\Theta (z)}=
\frac {2}{\pi i}\int \limits _{\mathbb R}
\bigg( \frac {1}{x-z}-\frac {x}{x^2+1}\bigg) d\nu (x).
$$
The condition (\ref {uslovie}) entails that
$$
\nu ({\mathbb R})<+\infty,
$$
and there is a limit $\lim \limits _{y\to +\infty}\Theta (iy)$; let us denote  it
by $\Theta (\infty)$. We have $|\Theta (\infty)|=1$
and $1-\overline{\Theta (\infty)}\Theta \in L^2({\mathbb R})$, with
$$
\|1-\overline {\Theta (\infty)}\Theta\|_{L^2({\mathbb R})}
=|1-\Theta (\infty)|\cdot\sqrt {\nu ({\mathbb R})}.
$$
The condition (\ref {uslovie}) is  equivalent to
$$
\int \limits _{\mathbb R}(1+|t|)^{q-2}d\nu (t)<\infty .
$$
The formula
$$
(Lf)(x)=\frac {1}{\sqrt \pi (x+i)}f\bigg(\frac {x-i}{x+i}\bigg)
$$
carries out the unitary mapping of the space $L^2({\mathbb T})$
to $L^2({\mathbb R})$ such that the Hardy space
$H^2({\mathbb D})$ transforms to the Hardy space $H^2({\mathbb C}_+)$.
Such a transformation turns 
the inclusion (\ref {inc})  into the relation
\begin{equation}
  \label{inc1}
  M_{1-\overline {\Theta (\infty)}\Theta } |_{\mathcal K} \in \mathfrak {S}_p,
\end{equation}
where $\mathcal{K} = H^2({\mathbb C}_+) \ominus e^{itz} H^2({\mathbb C}_+)$.
The Paley--Wiener space ${\mathcal PW}_{a}$
consists of all entire functions of the exponential type
at most $a$, the restriction of which to the real line
belongs to $L^2(\mathbb{R})$;
and, according to the classical Paley--Wiener theorem,
${\mathcal PW}_{a} =
e^{-iaz} H^2({\mathbb C}_+)\ominus e^{iaz}H^2({\mathbb C}_+)$.
In this case the inclusion  (\ref{inc1}) is  equivalent to the question, whether the transformation
of the Paley--Wiener space ${\mathcal PW}_{t/2}$ into the space
$L^2(\mathbb{R}, w(t)dt)$ on the real line with the weight  $w(t) =
|1-\overline {\Theta (\infty)}\Theta(t)|^2$ belongs to $\mathfrak {S}_p$.
This problem was solved in the paper \cite {Par}, with the following result obtained:
\smallskip

{\bf Theorem (O.G. Parfenov).} {\it For any
$p>0$ the embedding  operator $\mathcal J$
of the space ${\mathcal PW}_a,\ a>0$ into the space
$L^2(\mathbb{R}, w(t)dt)$ belongs to the class $\mathfrak {S}_p$ if and only if 
\begin{equation}
  \label{np}
\mathfrak{N}_p(w) = \sum\limits_k \bigg(\int \limits _k^{k+1}w(x)dx
\bigg)^{p/2}<\infty.
\end{equation}  }
The following estimate  follows immediately from the proof of the Parfenov theorem
(see also \cite {Anton}, where a similar result is obtained for the general
model spaces):
$$
\|{\mathcal J}\|_{\mathfrak {S}_p}^p\le \mathfrak{N}_p(w).
$$
\smallskip
\\
{\it \bf Step 4. Application of the Parfenov theorem.}
It follows from the embedding  $(1-\xi)^{-q}\in L^1(\mu)$ that the functional $\Phi$,
$$
\Phi(g) = \int \limits _{\mathbb T}
\bigg(\frac {1-\overline {\theta (1)}\theta (\xi)}{1-\xi}\bigg)^q g(\xi)d\xi,
\qquad g\in K_{\theta},
$$
is bounded on $K_{\theta}$, and
$|\Phi(g)| \le C(q) \mathfrak{M}_q(\mu) \|g\|_2$.
Note that for $q\in \mathbb{N}$ the value $\Phi(g)$
coincides with the radial limit $g^{(q-1)}(1)$
of the derivative of the order $q-1$ of the function $g$ at the point $z=1$.

Thus, the bounded functional $\Phi$ on $K_\theta$
is generated by the function $\big(\frac {1-\overline {\theta (1)}
\theta (\xi)}{1-\xi}\big)^q \in H^2(\mathbb{D})$. Strictly speaking,
the function $\big(\frac{1-\overline{\theta(1)}
\theta(\xi)}{1-\xi}\big)^q$ does not belong to the space
$K_\theta$, but one can easily show that the norm of its projection to the
subspace $\theta H^2$ is estimated via the norm of its projection
to $K_\theta$. Therefore,%, мы заключаем, что
\begin{equation}
  \label{otsenka}
\int\limits_{\mathbb{T}} \bigg|\frac{1-\overline{\theta(1)}\theta(\xi)}
{1-\xi}\bigg|^{2q} dm(\xi) \le \omega (\mathfrak{M}_q(\mu)),
\end{equation}
where $\omega (r) \to 0$ when $r \searrow 0$ (in fact, $\omega(r) \le C(q) r$,
but the explicit form of the function $\omega$ is not important for us).
Substituting the variable, we obtain
$$
\int \limits _{\mathbb R}|1-\overline {\Theta (\infty )}
\Theta (t)|^{2q}(|t|+1)^{2q-2}dt<\infty.
$$
Apply the H\"older inequality, we obtain
$$
\begin{aligned}
\int \limits _k^{k+1}|1-&\overline {\Theta (\infty)}\Theta (t)|^2dt\\
\le &  \bigg(\int \limits _k^{k+1}|1-\overline {\Theta (\infty)}\Theta (t)|^{2q}
(|t|+1)^{2q-2}dt\bigg)^{1/q}\bigg(\int \limits _k^{k+1}\frac {dt}{(|t|+1)^2}\bigg)^{1/q'} \\
\le & \frac {C^{1/q}}{(|k|+1)^{2/q'}}.
\end{aligned}
$$
Let us assume that $p>q'$; then
$$
\sum \limits _{k\in {\mathbb Z}}\bigg(\int \limits _k^{k+1}
|1-\overline {\Theta (\infty )}\Theta (t)|^2dt\bigg)^{p/2}
\le C^{p/2q}\sum \limits _{k\in {\mathbb Z}}\frac {1}{(|k|+1)^{p/q'}}<\infty.
$$
Thus, invoking the estimate (\ref {otsenka}) when $p>q',$ we obtain
$$
\sum \limits _{k\in \mathbb Z}\bigg(\int \limits _k^{k+1}
|1-\overline {\Theta (\infty )}\Theta (t)|^2dt\bigg)^{p/2}
\le \omega (\mathfrak{M}_q(\mu)),
$$
with some function $\omega$, $\omega (r)\searrow 0$ when $r\searrow 0$.
Then, applying the Parfenov theorem, we obtain the inclusion (\ref{inc1}).
Proposition 9 is proved completely.
\qed
\medskip

In the model of cocyclic perturbation considered here,
the relation $W_t-I\in\mathfrak S_p$ is equivalent to the inclusion
$\varphi _t(\tilde V)-\varphi _t(\tilde S)\in\mathfrak S_p$.
In conclusion to the section note that the difference
$\varphi _t(\tilde V)-\varphi _t(\tilde S)$ cannot belong to
the trace class $\mathfrak {S}_1$ for all $t\ge 0$ simultaneously.

\medskip

{\bf Proposition 10.} {\it For the class of  cocyclic perturbations described in Proposition~$5$,
the inclusion $\varphi _t(\tilde V)-\varphi _t(\tilde S)\in\mathfrak S_1$ 
for all $t\ge 0$ implies that $\theta$ is a unimodular constant.}

\medskip

{\bf Proof.}
It follows from Proposition 9, that the inclusion
$\varphi _t(\tilde V)-\varphi _t(\tilde S)\in \mathfrak {S}_1$
is  equivalent to $\mathfrak{N}_1(|1-
\overline{\Theta(\infty)}\Theta|^2) <\infty$ (see (\ref{np})).
It would result in
$$
\int\limits_{\mathbb{R}} |1-\overline{\Theta(\infty)}
\Theta(t)|\, dt \le \sum\limits_{k\in\mathbb{Z}}
\bigg(\int_{k}^{k+1}|1-\overline{\Theta(\infty)}
\Theta(t)|^2\,dt \bigg)^{1/2}< \infty,
$$
and therefore the function $1-\overline{\Theta(\infty)} \Theta$
should belong to the Hardy space $H^1$. But then
$\int_{\mathbb{R}} \big(1-\overline{\Theta(\infty)}
\Theta(t)\big)\, dt=0$, which is impossible since
${\rm Re}\, (1- \overline{\Theta(\infty)}  \Theta) >0$
almost everywhere on $\mathbb{R}$
for any nonconstant inner function $\Theta$. \qed

%%%%%%%%%%%%%%%%%%%%%%%%%%%%%%%%%%%%%%%%%%%%%%%%%%%%

\section {The case of an arbitrary spectral multiplicity}

Let $U$ be a unitary part in the Wold--Kolmogorov decomposition
(\ref {Wold1}) of the arbitrary cogenerator of the cocyclic perturbation.
Any unitary operator $U$ can be presented as an at most countable sum
$$
U=\oplus _kU_k,
$$
where the operators $U_k$ are unitarily equivalent to the operators of multiplication in the appropriate
spaces $L^2(\mu _k)$ and $\mu _k$ are measures on the circle $\mathbb T$,
$$
(U_kf)(z)=zf(z),\quad f\in L^2(\mu _k).
$$
Multiplying by positive weights, decreasing rapidly  near the point $1$,
we can choose measures $\mu _k$ such that the condition
\begin{equation}\label{multi}
\sum \limits _k\bigg(\int \limits _{\mathbb T}
\frac {d\mu _k(\xi)}{|1-\xi |^{q}}\bigg)^{1/q}<\infty
\end{equation}
holds for all $q>0$.
Let us define the inner functions $\theta _k$,
connected with the measures $\mu _k$ by the formula (\ref {Clark}). Condition (\ref {multi})
ensures that the product $\prod _k\theta _k$ converges to some inner function
$\theta $. Put
$$
\hat \theta _n=\prod _{k=1}^{n-1}\theta _k
$$
and define the cogenerator $\tilde V$ by the formula
$$
\tilde V=\tilde S+\sum \limits _n(\cdot ,\hat \theta _ng_n)
\hat \theta _n(1-\theta _n)-(\cdot ,\overline z)(1-\overline {\theta (1)}\theta ),
$$
where
$$
g_n(z)=\frac {\theta _n(z)-\theta (0)}{z(1-\theta _n(0))}.
$$
\medskip

{\bf Proof of Theorem 1.}
The operator $V=\tilde V|_{K}$ is diagonal with respect to
the orthogonal decomposition
$K=\oplus _k\hat \theta _kK_{\theta _k}\oplus \theta K$.
Condition (\ref {multi}) and Proposition 8 entail that
$$
\varphi _t(V)-\varphi _t(S)\in \mathfrak {S}_1,\qquad t\ge 0.
$$
The same condition (\ref {multi}) and Proposition 9 provide the inclusion
$$
\varphi _t(\tilde V)-\varphi _t(\tilde S)\in \mathfrak {S}_p,\qquad t\ge 0,
$$
for $p>q'$. Since the condition (\ref {multi}) 
holds for arbitrarily large values of $q$ by the choice 
of measures, we have $\varphi _t(\tilde V)-\varphi _t(\tilde S)\in 
\mathfrak {S}_p$ for any $p>1$.
\qed
\medskip

{\bf Proof of Theorem 2.}
Let $U$ be a cogenerator of an arbitrary semigroup
of unitary operators, being a unitary part in the Wold--Kolmogorov decomposition
of the cocyclic perturbation. Then, there is
an operator $\Delta$, belonging to
all classes $\mathfrak {S}_p$ for $p>1$, that  the
perturbation $U+\Delta$ has a singular spectrum (see \cite {Kato}).
Moreover,
$$
\varphi _t(U+\Delta )-\varphi _t(U)\in \mathfrak {S}_p,\ t\ge 0.
$$
A detailed proof of the last statement is given
in \cite{ABK} (proof of Theorem 1.3).
To complete the proof it is sufficient to apply Theorem 1.
\qed

\end{document}